\documentclass[11pt,a4paper]{amsart}

\usepackage[toc,page]{appendix}
\usepackage{fancyhdr}
\usepackage[hypertex]{hyperref} 

\usepackage{amsmath,amsthm, amscd, amssymb, amsfonts}
\usepackage{epsfig}
\usepackage{amsmath,amsthm,amssymb, amscd,enumerate}
\usepackage{pb-diagram}

\newcommand{\C}{{\mathbb C}}
\numberwithin{equation}{section}

\theoremstyle{plain}
\newtheorem{theorem} {Theorem} [section]
\newtheorem*{theoremSN} {Theorem}
\newtheorem{lemma} [theorem] {Lemma}
\newtheorem{corollary} [theorem] {Corollary}

\theoremstyle{definition}
\newtheorem{definition} [theorem] {Definition}
\newtheorem{remark} [theorem] {Remark}

\newtheorem*{exampleSN} {Example}

\renewcommand \parallel {/\kern-3pt/}

\newcommand \g {\mathfrak{g}}
\newcommand \z {\mathfrak{z}}
\newcommand \n {\mathfrak{n}}
\newcommand \gl {\mathfrak{gl}}
\newcommand \h {\mathfrak{h}}

\renewcommand \k {\mathrm{k}}

\newcommand \End {\operatorname{End}}

\newcommand \im {\operatorname{im}}

\begin{document}

\title[Faithful representations of current Heisenberg Lie algebras]
{Faithful representations of minimal dimension of current
Heisenberg Lie algebras }
\author{Leandro Cagliero}
\address{CIEM-FaMAF, Universidad Nacional de C\'ordoba}
\email{cagliero@famaf.unc.edu.ar}

\author{Nadina Rojas}
\address{CIEM-FaMAF, Universidad Nacional de C\'ordoba}
\email{nrojas@famaf.unc.edu.ar}

\begin{abstract}
Given a Lie algebra $\g$ over a field of characteristic zero $\k$,
let $\mu(\g)=\min\{\dim \pi: \pi\text{ is a faithful representation of }\g\}$.
Let $\h_{m}$ be
the Heisenberg Lie algebra of dimension $2m+1$ over $\k$ and let
$\k[t]$ be the polynomial algebra in one variable. Given
$m\in\mathbb{N}$ and $p\in\k[t]$, let $\h_{m,p}=\h_m\otimes
\k[t]/(p)$ be the current Lie algebra associated to $\h_m$ and
$\k[t]/(p)$, where $(p)$ is the principal ideal in $\k[t]$ generated by $p$.
In this paper we prove that $ \mu(\h_{m,p}) = m \deg p + \left
\lceil 2\sqrt{\deg p} \right\rceil. $
\end{abstract}

\date{\today}

\maketitle





\section{Introduction}


Let $\k$ be a fixed field of characteristic zero.
In this paper, all Lie algebras, associative algebras,
Hom and tensors  are considered over $\k$,
unless otherwise is explicitly mentioned.

By Ado's theorem every
finite dimensional Lie algebra has a finite dimensional faithful representation (see for instance \cite{JN}).
However, given a Lie algebra $\g$, it is in general very difficult to compute the number
\[
\mu(\g)=\min\{\dim V: (\pi,V)\text{ is a faithful representation of }\g\}.
\]

The problem of computing the value of $\mu(\g)$, or bounds for it, gained interest since
Milnor \cite{Mi} posed the question of which are the finite groups that occur as
fundamental groups of complete affinely flat manifolds;
in particular whether they are the polycyclic-by-finite groups.
There are many articles giving an affirmative answer to Milnor's question
under some additional hypothesis, see for instance \cite{Au}, \cite{De}, \cite{GK}, \cite{GM}, etc.
However, the answer to the original Milnor's question is negative in both directions.
On the one hand, Margulis \cite{Ma} gave the first complete affinely flat manifold whose
fundamental group do not have a polycyclic subgroup of finite index.
On the other hand, Benoist \cite{BY} and Burde and Grunewald \cite{BG}
found the first examples of nilpotent Lie groups without
any left-invariant affine structures.
These examples are achieved by finding nilpotent Lie algebras $\g$ such that $\mu(\g)>\dim(\g)+1$.

Very little is known about the function $\mu$.
In particular the value of $\mu$ is known only for a few families of Lie algebras,
among them, reductive Lie algebras over algebraically closed fields (see \cite{BW}),
abelian Lie algebras, and Heisenberg Lie algebras.
A brief account of some known results is the following:
\begin{enumerate}[1.]
  \item If $\g$ is abelian, then $\mu(\g) = \left \lceil 2 \sqrt{\dim(\g) - 1}\right\rceil$.
  Here $\lceil a \rceil$ is the closest integer that is greater than or equal to $a$.
  This result is due to Schur \cite{S}, for $\k=\mathbb{C}$, and to Jacobson \cite{J2} for arbitrary $\k$
  (see also \cite{MM}).
\item If
      $\h_m$ is the Heisenberg Lie algebra of dimension $2m+1$, then $\mu(\h_m) = m + 2$, (see \cite{BU2}).
\item If
      $\g$ is a $k$-step nilpotent Lie algebra, then
      $\mu(\g) \leq 1 + (\dim \g)^k$. This is part of Birkhoff´s embedding theorem
      (see \cite{Bi} and also \cite{REE}).
      In addition, if  $\g$ is $\mathbb{Z}$-graded, then $\mu(\g) \le \dim \g$ (see \cite{BU2}).
\item If
      $\g$ is a filiform Lie algebra then $\mu(\g) \ge \dim \g$ (see \cite{BY}).
       The equality holds if  $\dim \g < 10$ (see \cite{BU1}).
\end{enumerate}

In this paper we compute the value of $\mu$ for a whole family of
current Lie algebras associated to Heisenberg Lie algebras.
Let $\k[t]$ be the polynomial algebra in one variable and
given  $p\in\k[t]$, let $(p)$ denote the principal ideal generated by $p$ in $\k[t]$.
Given $m\in\mathbb{N}$ and a non-zero polynomial $p\in\k[t]$ let
\[
\h_{m,p}=\h_m\otimes\k[t]/(p)
\]
be the Lie algebra over $\k$ with bracket
$[X_1\otimes q_1,X_2\otimes q_2]=[X_1,X_2]\otimes q_1q_2$,
$X_i\in\h_m$, $q_i\in\k[t]/(p)$, $i=1,2$.
It is clear that
\[
\dim\h_{m,p}=(2m+1)\deg p.
\]
The set of Lie algebras $\h_{m,p}$ constitute a family of 2-step nilpotent Lie algebras that contains
the following two subfamilies.

\medskip

\noindent \emph{Truncated Heisenberg Lie algebras.}
This is the subfamily corresponding to $p=t^k$,
$k\in\mathbb{N}$.
These Lie algebras appear in the literature associated to
the {Strong Macdonald Conjetures} \cite{Mac}.
Some articles dealing with these conjectures are \cite{FGT}, \cite{HW}, \cite{Ku}, \cite{T}.

\medskip

\noindent \emph{Heisenberg Lie algebras over finite extensions of $\k$.}
This is the subfamily corresponding to polynomials $p$ that are irreducible over $\k$.
In this case $\h_{m,p}$ is the Lie algebra obtained by restricting scalars to $\k$ in
the Heisenberg Lie algebra over the field $K_p=\k[t] /(p)$.

The main result of this paper is the following theorem.

\begin{theoremSN} \label{mudeHm,p}
Let $m\in\mathbb{N}$ and $p\in\k[t]$, $p\ne0$. Then
\[
\mu(\h_{m,p}) = m \deg p + \left \lceil 2\sqrt{\deg p} \right\rceil.
\]
\end{theoremSN}

In order to prove the inequality
$\mu(\h_{m,p}) \ge m \deg p + \left \lceil 2\sqrt{\deg p} \right\rceil$,
we first prove Theorem \ref{corleandro} in which we obtain
some fine information about the structure of a faithful representation of an abelian Lie algebra.
From Theorem \ref{corleandro} it is straightforward to obtain
the lower bound part of the Theorem of Schur mentioned above (see Corollary \ref{Coro.Schur}).
One might expect that this result will help to obtain lower bounds of $\mu(\g)$ for other families
of nilpotent Lie algebras.

On the other hand, the proof of
$\mu(\h_{m,p}) \le m \deg p + \left \lceil 2\sqrt{\deg p} \right\rceil$
is done by explicitly constructing faithful
representations of $\h_{m,p}$ of minimal dimension.

\begin{exampleSN}
Let $\h_{m}(\C)$ be the Heisenberg Lie algebra of dimension
$2m+1$ over the complex numbers
and let
$\h_{m}(\C)_{\mathbb{R}}$ be $\h_{m}(\C)$ viewed as a Lie algebra over $\mathbb{R}$.
We know that $\mu(\h_{m}(\C))=m+2$ and
the faithful representation of $\h_{m}(\C)$ in $\C^{m+2}$
yields a faithful representation of $\h_{m}(\C)_{\mathbb{R}}$
of dimension $2m+4$.
However, from the above theorem we obtain that $\mu(\h_{m}(\C))_{\mathbb{R}}=2m+3$.
If $\{X,Y,Z\}$ is the basis of $\h_{1}(\C)$, with $[X,Y]=Z$,
then
\[
\pi\big((x_1+ix_2)X+(y_1+iy_2)Y+(z_1+iz_2)Z\big)
 =\left(
            \begin{smallmatrix}
                   0 & x_1 & x_2 & z_1  & z_2 \\[1mm]
                   0 & 0   & 0   & y_1  & y_2 \\[1mm]
                   0 & 0   & 0   & -y_2 & y_1 \\[1mm]
                   0 & 0   & 0   & 0    & 0   \\[1mm]
                   0 & 0   & 0   & 0    & 0
            \end{smallmatrix}\right)
\]
is a faithful representation of $\h_{1}(\C)_{\mathbb{R}}$ in $\mathbb{R}^5$.
\end{exampleSN}


\section{Preliminaries}


Given a finite dimensional vector space $V$, a representation of
an associative algebra $A$ on $V$ is an associative algebra homomorphism
$\pi:A\to\text{End}(V)$ and similarly, a representation of a Lie
algebra $\g$ on $V$ is Lie algebra homomorphism $\pi:\g\to\gl(V)$.
A \emph{nil-representation} is, by definition, a representation
whose image is contained in the set of nilpotent endomorphisms.
All representations $(\pi,V)$ considered in this paper will be
finite dimensional. Thus, given a basis of $V$, we can express
each operator $\pi(X)$ by its associated matrix. When the basis is
fixed we shall denote this matrix also by $\pi(X)$. The space of
rectangular matrices of size $m\times n$ with entries in $\k$ will
be denoted by $M_{m, n}(\k)$.

A representation is \emph{faithful} if it is an injective homomorphism.
Considering the left action, it is clear that any associative algebra with unit
of dimension $n$ has a faithful representation of dimension $n$.
On the other hand, a theorem due to Ado states that
any finite dimensional Lie algebra has a  finite dimensional faithful representation
(see for instance \cite{JN}).
Given a Lie algebra $\g$, let
\[
\mu(\g)=\min\{\dim V: (\pi,V)\text{ is a faithful representation of }\g\}.
\]

Given a Lie algebra $\g$ and a commutative associative algebra $A$
the tensor product $\g \otimes A$ has a Lie algebra structure with bracket
\[
[X_1\otimes a_1,X_2\otimes a_2]=[X_1,X_2]\otimes a_1a_2,
\]
$X_i\in\g$, $a_i\in A$, $i=1,2$. This Lie
algebra is known as the \emph{current Lie algebra} associated to
$\g$ and $A$ (see for instance \cite{GR}, \cite{Z2}).

\begin{remark}\label{rmk. field extension}
Note that $\g \otimes A$ could
be viewed as a Lie algebra over the algebra $A$ but, as we have already mentioned,
we look at it as a Lie algebra over $\k$.
Occasionally we shall make an exception when
$A=K$ is a field extension of $\k$.
In this case, we shall denote by $\g_K$ the Lie algebra
$\g \otimes K$ whenever it is viewed as a Lie algebra over $K$.
\end{remark}

If
$(\pi, V)$ is a representation of a Lie algebra
$\g$ and
$(\rho , W)$ is a representation of a commutative associative algebra
$A$, then it is clear that
$\pi \otimes \rho : \g \otimes A \rightarrow \mathfrak{gl}(V \otimes W)$, given by
$(\pi \otimes \rho)(X \otimes a) = \pi(X) \otimes \rho(a)$, is a representation of the Lie algebra
$\g \otimes A$.
It is clear that if $\pi$ and $\rho$ are injective, then so it is $\pi \otimes \rho$.
Therefore, if $A$ has an identity, by considering the regular representation of $A$ we obtain
that
\begin{equation}\label{eq.first inequality}
\mu(\g \otimes A) \leq \mu(\g)\dim A.
\end{equation}
If $A=K$ is a field extension of $\k$
and we look at $\g_K=\g\otimes K$ as a Lie algebra over $K$,
then the above
construction also yields a $K$-representation of $\g_K$
and in this case one has
\begin{equation}\label{eq.first inequality for filed extension}
\mu(\g_K) \leq \mu(\g).
\end{equation}

The Heisenberg Lie algebra $\h_m$ is the $\k$-vector space of dimension $2m+1$ with a basis
 $\{X_1,\ldots ,X_m, Y_1,\ldots ,Y_m,Z\}$ whose only non-zero brackets are
 \[
 [X_i, Y_i]= Z.
 \]
It is clear that the center of $\h_m$ is $\mathfrak{z}(\h_m)=\k Z$.
This Lie algebra has a well known faithful representation $(\pi_0,\k^{m+2})$
that, in terms of the canonical basis of $\k^{m+2}$,
is given by
\begin{equation}\label{Rep of heisenberg}
\textstyle
\pi_0\big(\sum_{i=1}^m x_i X_i + \sum_{i=1}^m y_i Y_i + z Z \big) =
      \left(\begin{smallmatrix}
           0 & x_1 & \ldots & x_m & z  \\[.8mm]
              &     &        &     & y_1 \\
              &     & 0      &     &  \vdots  \\[.8mm]
              &     &        &       & y_m  \\[.8mm]
              &     &        &       & 0
      \end{smallmatrix}\right),\quad x_i,y_i,z\in\k.
\end{equation}
It is known that any other faithful representation of $\h_m$ has dimension greater than or equal to $m+2$
(see \cite{BU2}) and thus $\mu(\h_m)=m+2$.

Let $\k[t]$ be the polynomial algebra in one variable,
let $p=a_0+\dots+a_{d-1}t^{d-1}+t^d$ be a non-zero monic polynomial and
let $(p)$ be the principal ideal generated by $p$.
The regular representation $\rho$ of the quotient algebra $\k[t]/(p)$
is expressed, in terms of the canonical basis $\{1,t,\dots,t^{d-1}\}$, by
$\rho(t^i) = P^i$, where
\begin{equation}\label{matriz P0}
P =
       \left(\begin{smallmatrix}
       \;0\;   &      &        &          & -a_0    \\[1mm]
       \;1\;   & \;0    &        &        & -a_1    \\[-2mm]
           & \;1\;    & \ddots &          &  \cdot       \\[-.5mm]
           &      &          & \;         &   \cdot       \\[-4mm]
           &      & {\ddots} & \;0        &   \cdot       \\[2mm]
           &      &        & \;1\;        & -a_{d-1}
      \end{smallmatrix}\right)
\end{equation}
is the matrix associated to $p$.

For any $m\in\mathbb{N}$ and any non-zero polynomial $p\in\k[t]$ let
\[
\h_{m,p}=\h_m\otimes\k[t]/(p)
\]
be the \emph{current Heisenberg Lie algebra} associated to $m$ and $p$.
Note that  \eqref{eq.first inequality} yields
\[
\mu(\h_{m,p})\le (m+2)\deg p=m\deg p+2\deg p.
\]
The main goal of this paper is to  prove that in fact
\[
\mu(\h_{m,p})=m\deg p+ \left \lceil 2\sqrt{\deg p} \right\rceil.
\]
Note that this extends the known result $\mu(\h_m)=m+2$.

\medskip

We now recall and prove some results, needed in the following sections,
about finite dimensional representations of nilpotent Lie algebras.

Let $\n$ be finite dimensional nilpotent Lie algebra and let
  $(\pi,V)$ be a finite dimensional representation of $\n$.
  If $\k$ is algebraically closed,
  a well known theorem of Zassenhaus (see \cite{JN}) states that
$V$ can be decomposed as
$$V = V_1 \oplus V_2 \oplus \dots \oplus V_r,$$
so that, for all $X\in\n$ and $i=1,\dots,r$,
$\pi(X)|_{V_i}$ is an scalar $\lambda_i(X)$ plus a nilpotent operator
$N_i(X)$ on $V_i$.

A consequence of this result is
is that every representation of a nilpotent Lie algebra has
a Jordan decomposition.

\begin{theorem}
  Let $\n$ be finite dimensional nilpotent Lie algebra and let
  $(\pi,V)$ be a finite dimensional representation of $\n$.
  For each $X\in\n$ let $\pi_S(X)$ and $\pi_N(X)$ be, respectively,
  the semisimple and nilpotent part of the additive Jordan decomposition of $\pi(X)$.
  Then $(\pi_S , V)$ and $(\pi_N , V)$ are representations of $\n$.
  Moreover, $\pi_S(X)=0$ for all $X\in[\n,\n]$.
\end{theorem}

\begin{proof}
We first assume that $\k$ is algebraically closed.
By Zassenhaus' theorem
$$V = V_1 \oplus V_2 \oplus \dots \oplus V_r,$$
and $\pi(X)|_{V_i}$ is an scalar $\lambda_i(X)$ plus a nilpotent operator
$N_i(X)$ on $V_i$,
for all $X\in\n$ and $i=1,\dots,r$.
Since $\lambda_i(X)=\text{trace}(\pi(X)|_{V_i})/\dim(V_i)$ it follows that
$\lambda_i$ is linear and $\lambda_i(X)=0$ for all $X\in[\n,\n]$.
Hence
$\lambda_i$ is a Lie algebra homomorphism.
Moreover, since $N_i=\pi|_{V_i}-\lambda_i$ and $\lambda_i(X)$ is a scalar,
it follows that $N_i$ is a Lie algebra homomorphism.
Since $\pi_S|_{V_i}=\lambda_i$ and $\pi_N|_{V_i}=N_i$ the theorem follows for
$\k$ algebraically closed.

Now assume $\k$ is arbitrary (with char$(\k)=0$).
Let $\bar\k$ be the algebraic closure of $\k$,
$\bar\g=\g\otimes_\k\bar\k$,
$\bar V=V\otimes_\k\bar\k$ and $\bar\pi=\pi\otimes_\k 1$.
We know that
$\bar \pi_S$ and $\bar \pi_N$ are Lie algebras homomorphisms and
we must prove that $\pi_S$ and $\pi_N$ are Lie algebras homomorphisms.

Let $B = \{v_1, \dots, v_n\}$ be a basis of $V$,
$\bar{B} = \{v_1 \otimes 1, \dots, v_n \otimes 1\}$
the corresponding basis of $\bar{V}$ and,
given an operator $T$ on $V$ (resp. on $\bar{V}$) let
$[T]_B$ (resp. $[T]_{\bar B}$) be the associated matrix with respect to the basis $B$ (resp. $\bar B$).
Since $[\pi(X)]_B= [\bar{\pi}(X\otimes_\k 1)]_{\bar{B}}$ for all $X \in \g$,
it follows that
\begin{equation}\label{eq.jordan decomp}
[\bar{\pi}(X\otimes_\k 1)]_{\bar{B}} = [{{\pi}_S}({X})]_{{B}} + [{{\pi}_N}({X})]_{{B}}.
\end{equation}
Since $[{{\pi}_S}({X})]_{{B}}$ and $[{{\pi}_N}({X})]_{{B}}$
are respectively semisimple and  nilpotent matrices that commute, it follows that
\eqref{eq.jordan decomp} is the Jordan decomposition of the matrix
$[\bar{\pi}(X\otimes_\k 1)]_{\bar{B}}$, that is
$[\pi_S(X)]_{{B}}=[\bar{\pi}_S(X\otimes_\k 1)]_{\bar{B}}$ and
$[\pi_N(X)]_{{B}}=[\bar{\pi}_N(X\otimes_\k 1)]_{\bar{B}}$.
Therefore $\pi_S$ and $\pi_N$ are Lie algebras homomorphisms.
\end{proof}

\begin{definition}
  Let $(\pi,V)$ be a finite dimensional representation of a
  finite dimensional nilpotent Lie algebra $\n$.
  We  call the representations $(\pi_S , V)$ and $(\pi_N , V)$
  the \emph{semisimple part }and the \emph{nilpotent part} of $(\pi,V)$ respectively.
\end{definition}

For certain nilpotent Lie algebras it is enough to consider nilrepresentations
in the definition of $\mu$ as the following theorem shows.

\begin{theorem}\label{thm.nilrep}
  Let $\n$ be a finite dimensional nilpotent Lie algebra such that
  the center $\mathfrak{z}(\n)$ is contained in $[\n,\n]$,
  and let $(\pi,V)$ be a finite dimensional representation of $\n$.
  Then $(\pi,V)$ is faithful if and only if $(\pi_N , V)$ is faithful.
\end{theorem}

\begin{proof}
We only need to show that if $\pi$ is injective, then so it is $\pi_N$.
Let us assume that $\pi$ is injective and let $X_0\in\n$ such that $\pi_N(X_0)=0$.
Since $\pi_S|_{[\n,\n]} = 0$ it follows that
\begin{align*}
 \pi([X_0 , X]) &= \pi_N([X_0 , X])+ \pi_S([X_0 , X])\\
                      &= [\pi_N(X_0) , \pi_N(X)]  \\
                      &= 0
\end{align*}
for all  $X \in \n$.
Since $\pi$ is injective we obtain that $X_0 \in \mathfrak{z}(\n)\subset [\n,\n]$.
Finally, since $\pi|_{[\n,\n]}=\pi_N|_{[\n,\n]}$it follows that  $\pi(X_0) = 0$
and therefore $X_0 = 0$.
\end{proof}


\section{A family of  representations of $\h_{m,p}$ \\ and the upper bound for $\mu(\h_{m,p})$}


Let $m\in\mathbb{N}$ be a fixed natural number and let $p=a_0+\dots+a_{d-1}t^{d-1}+t^d$
be a fixed non-zero monic polynomial of degree $d$.
If
$(\pi_0, \k^{m+2})$ is the faithful representation of $\h_m$ defined in \eqref{Rep of heisenberg}
and
$(\rho , \k[t]/(p))$ is the regular representation of  $\k[t]/(p)$,
then the representation
\[
\pi_0 \otimes \rho : \h_{m,p} \rightarrow \mathfrak{gl}(\k^{m+2} \otimes \k[t]/(p))
\]
is expressed, in terms of the canonical basis of the tensor product $\k^{m+2} \otimes \k[t]/(p)$, by
the blocked matrix
\begin{multline*}
(\pi_0 \otimes \rho )
\left(
\sum_{i=1}^m  X_i\otimes q_{1,i}(t) +
\sum_{i=1}^m Y_i\otimes q_{2,i}(t) +
Z\otimes q_3(t)
\right) \\
=
\left(\begin{smallmatrix}
           0 & q_{1,1}(P) & \dots & q_{1,m}(P) & q_{3}(P) \\[1mm]
             &                     &        &                     & q_{2,1}(P) \\[1mm]
             &                     &    0   &                     &  \vdots    \\[1mm]
             &                     &        &                     & q_{2,m}(P) \\[1mm]
             &                     &        &                     & 0
      \end{smallmatrix}\right)
\end{multline*}
of size $(m+2)d$.
We shall now construct a family of representations of $\h_{m,p}$ that contains $\pi_0 \otimes \rho$.

\begin{definition}
Given two natural numbers $a$ and $b$ and two matrices
$A \in M_{a,d}(\k)$ and $B \in M_{d ,b}(\k)$,
let $(\pi_{A,B},\k^{md+a+b})$ be the representation of
$\h_{m,p}$ that, in terms of the canonical basis of $\k^{md+a+b}$,  is given by the blocked matrix
\begin{multline*}
\pi_{A,B}
\left(
\sum_{i=1}^m  X_i\otimes q_{1,i}(t) +
\sum_{i=1}^m Y_i\otimes q_{2,i}(t) +
Z\otimes q_3(t)
\right) \\
=
\left(\begin{smallmatrix}
           0\; & \;A\;q_{1,1}(P)\; & \;\dots\; & \;A\;q_{1,m}(P)\; & \;A\;q_{3}(P)\,B \\[2mm]
             &                     &        &                     & \;q_{2,1}(P)\,B \\[2mm]
             &                     &    0   &                     &  \;\vdots    \\[2mm]
             &                     &        &                     & \;q_{2,m}(P)\,B \\[2mm]
             &                     &        &                     & 0
      \end{smallmatrix}\right),
\end{multline*}
whose block structure is depicted below:
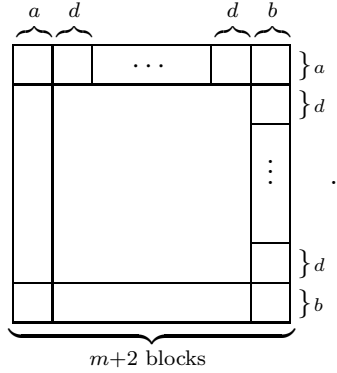
\begin{center}
\setlength{\unitlength}{15pt}
\begin{picture}(9,10)(0,-1.8)
\put(1.9,-.1){$\underbrace{\rule{108pt}{0pt}}_{m + 2 \text{ blocks }}$}
\put(2,0){\line(0,1){7}}
\put(3,0){\line(0,1){7}}
\put(9,0){\line(0,1){7}}
\put(8,0){\line(0,1){7}}
\put(2,0){\line(1,0){7}}
\put(2,1){\line(1,0){7}}
\put(2,6){\line(1,0){7}}
\put(2,7){\line(1,0){7}}
\put(4,6){\line(0,1){1}}
\put(5,6.3){$\cdots$}
\put(7,6){\line(0,1){1}}
\put(2.1,7.1){$\tiny \overbrace{}^{a}$}
\put(3.1,7.1){$\tiny \overbrace{}^{d}$}
\put(7.1,7.1){$\tiny \overbrace{}^{d}$}
\put(8.1,7.1){$\tiny \overbrace{}^{b}$}
\put(8,2){\line(1,0){1}}
\put(8.4,3.5){$\vdots$}
\put(8,5){\line(1,0){1}}
\put(9.1,6.3){$\left.\rule{0pt}{0pt}\right\}$}
\put(9.6,6.3){\tiny $a$}
\put(9.1,5.3){$\left.\rule{0pt}{0pt}\right\}$}
\put(9.6,5.3){\tiny $d$}
\put(9.1,1.3){$\left.\rule{0pt}{0pt}\right\}$}
\put(9.6,1.3){\tiny $d$}
\put(10,3.5){.}
\put(9.1,0.3){$\left.\rule{0pt}{0pt}\right\}$}
\put(9.6,0.3){\tiny $b$}
\end{picture}
\end{center}
\end{definition}
The following facts are not difficult to prove:
\begin{enumerate}[(1)]
\item $(\pi_{A,B},\k^{md+a+b})$ is a representation of $\h_{m,p}$ for any $A$ and $B$.
\medskip
\item If $a=b=d$ and $A$ and $B$ are the identity matrix, then $\pi_{A,B}=\pi_0 \otimes \rho$.
\medskip
\item Let $\k[P]\subset M_{d,d}(\k)$ be the subalgebra generated by $P$ in $M_{d,d}(\k)$
and  let $\beta_{A,B}:\k[P]\to M_{a,b}(\k)$ be the linear map given by
$\beta_{A,B}(q(P))=Aq(P)B$.
Then $\pi_{A,B}$ is injective if and only if $\beta_{A,B}$ is injective.
\end{enumerate}

Now the problem of finding the minimal possible dimension for a faithful representation
among the family $\pi_{A,B}$ reduces to finding the minimal value of $a+b$
among the pairs $(a,b)$
for which there exists two matrices
$A \in M_{a,d}(\k)$ and $B \in M_{d ,b}(\k)$
such that $\beta_{A,B}$ is injective.
Since $\dim(\k[P])=d$ and $\dim(M_{a,b}(\k))=ab$
it follows necessarily that $ab\ge d$.
We shall see now
that this condition is sufficient as well.
Let
$A \in M_{a,d}$ and
$B \in M_{d,b}$ be the following matrices
\begin{equation}\label{eq:matrices A y B}
 A_{ij} =
      \begin{cases}
      1, & \text{ if } j = d-(a-i)b\,;\\
      0, & \text{ otherwise; }
      \end{cases}
      \qquad
B_{ij} =
      \begin{cases}
      1, & \text{ if } i = j\,; \\
      0, & \text{ otherwise. }
      \end{cases}
\end{equation}
That is $B =
           \left(\!\!\!
      \begin{smallmatrix}       \\[4mm]   \\[4mm] \end{smallmatrix}
           \right.
      \underbrace{
      \begin{smallmatrix}
        1 &        &       \\
          & \ddots &     \\
          &        & 1
      \end{smallmatrix}}_d\;
      \underbrace{
       \begin{smallmatrix}
        0   & \cdots & 0     \\
     \vdots & \ddots & \vdots    \\
        0   & \cdots & 0
      \end{smallmatrix}}_{b-d}
           \left.
      \begin{smallmatrix}       \\[4mm]   \\[4mm]
      \end{smallmatrix}
           \!\!\!\right)$,
if $b\ge d$;
$B =
\left(
      \begin{smallmatrix}
        1 &        &       \\
          & \ddots &     \\
          &        & 1   \\[1mm]
        0   & \cdots & 0     \\
     \vdots & \ddots & \vdots    \\
        0   & \cdots & 0
      \end{smallmatrix}\right)\!\!
            \begin{array}{l}
        \left.\rule[-5mm]{1pt}{0pt}\right\}{\scriptstyle{b}}      \\
        \left.\rule[-5mm]{1pt}{0pt}\right\}{\scriptstyle{d-b}}      \\
      \end{array}
     $
if $d\ge b$ and,
for instance,
if
$d = 6$, $a = 4$ and
$b = 2$ we have
$A =
\left(
\begin {smallmatrix}
       0&0&0&0&0&0\\[.6mm]
       0&1&0&0&0&0\\[.6mm]
       0&0&0&1&0&0\\[.6mm]
       0&0&0&0&0&1\end {smallmatrix}\right)
$.

\begin{theorem}\label{Pro:B(AB)}
  If $ab\ge d$ and $A$ and $B$ are the matrices defined above, then
  $\beta_{A,B}$ is injective.
\end{theorem}
\begin{proof}
Let $q \in \k[t]$ be a monic polynomial such that $\deg(q)<d$ and let us show that
$Aq(P)B \ne 0$.
It is easy to prove by induction that,
for $1 \leq j \leq d - k$, one has
\begin{equation*}
(P^k)_{ij} =
       \begin{cases}
      1, & \text{ if } i = j + k;\\
      0, & \text{ if } i \ne j + k.
      \end{cases}
\end{equation*}
This implies that, for $1 \leq j \leq d - k$, one has
\begin{equation*}
(P^k B)_{ij} =
\begin{cases}
       1, & \text{ if } i = j + k;\\
       0, & \text{ if } i \ne j + k;
\end{cases}
\end{equation*}
and therefore, for  $1 \leq j \leq d - k$, one has
\begin{equation*}
(AP^k B)_{ij} =
\begin{cases}
      1 & \text{ if } j+k = d - (a - i)b;\\
      0 & \text{ if } j+k \neq d - (a - i)b.
      \end{cases}
\end{equation*}

Let $k_0 = \deg(q)$ and let us assume for a moment that for
there exist $i_0$ and $j_0$ such that
\begin{enumerate}[1.]
\item $1 \leq j_0 \leq \min\{b , d-k_0\}$,
\item $1\leq i_0 \leq a$ and
\item $j_0 + k_0 = d - (a - i_0)b$.
\end{enumerate}
Now we have
$$(AP^k B)_{i_0,j_0} =
\begin{cases}
      1, & \text{ if } k = k_0;\\
      0, & \text{ if } k = 0, \dots, k_0-1;
      \end{cases}$$
and this implies that
$(A q(P) B)_{i_0,j_0} =1$ and thus $A q(P) B \ne0$ as we wanted to prove.

In order to prove the
existence of $i_0$ and $j_0$ satisfying properties 1--3 above,
we shall see that
$t_0 =\left\lceil \frac{d-k_0}{b} - 1 \right\rceil$
satisfies $0\le t_0\le a-1$ and
$1 \leq d - k_0 - t_0 b \leq b$.
Once this is proved, $i_0 = a - t_0$ and
$j_0 = d - k_0 - t_0 b$ have the desired properties.

It is clear that $t_0\ge0$.
Since
$k_0<d \leq ab$  it follows that
$\left\lfloor \frac{d - k_0 -1}{b}\right\rfloor \leq a-1$.
It is easy to see that $\left\lceil \frac{x}{y} - 1 \right\rceil \leq \left\lfloor \frac{x - 1}{y} \right\rfloor$
for all
$x, y \in \mathbb{N}$, in particular
\[
t_0=
\left\lceil \frac{d-k_0}{b} - 1\right\rceil \leq
\left\lfloor \frac{d - k_0 - 1}{b} \right\rfloor\le a-1.
\]
Additionally
$$\frac{d-k_0}{b} - 1 \leq t_0 \leq \frac{d - k_0 - 1}{b}$$
and thus
$b \geq d - k_0 - t_0 b \geq 1$.
\end{proof}

\begin{corollary}
Let
$p \in \k[t]$ be non-zero polynomial. Then, for all
$m \in \mathbb{N}$
$$\mu(\h_{m,p}) \leq m \deg(p) + \left\lceil 2 \sqrt{\deg(p)} \right\rceil.$$
\end{corollary}

\begin{proof}
By Proposition \ref{Pro:B(AB)},
$\mu(\h_{m,p}) \leq m \deg(p) + a + b$ for all $a$ and $b$ such that $ab>\deg(p)$.
Since
\begin{equation}\label{eq:parte entera}
\min\{a+b:a,b\in \mathbb{N} \text{ and } ab\ge d\}= \left \lceil 2\sqrt{d} \right\rceil
\end{equation}
 for all $d\in \mathbb{N}$,
 we obtain the desired inequality.
\end{proof}


\section{The lower bound for $\mu(\h_{m,p})$}


From the inequality \eqref{eq.first inequality for filed extension} we know that
\[
\mu\Big(\big(\g \otimes \k[t]/(p)\big)_K\Big)\le\mu\big(\g \otimes \k[t]/(p)\big).
\]
Since
$
\big(\g \otimes \k[t]/(p)\big)_K \simeq \g_K\otimes_K K[t]/(p)
$
as Lie algebras over $K$, we obtain
\[
\mu\big(\g_K\otimes_K K[t]/(p)\big)\le\mu\big(\g \otimes \k[t]/(p)\big).
\]
Therefore, in order to obtain a lower bound for $\mu(\h_{m,p})$ we may assume that
$\k$ is algebraically closed.
In this case
$p = (t - b_1)^{d_1} \dots (t - b_r)^{d_r}$ for different $b_l \in \k$
and
\begin{align*}
\h_{m,p}
&\simeq \bigoplus_{l = 1}^r \h_m \otimes \k[t]/\big((t - b_l)^{d_l}\big) \\
&\simeq \bigoplus_{l = 1}^r \h_{m,t^{d_l}}.
\end{align*}
This section is devoted to prove the following theorem.
\begin{theorem}\label{TeoPrincipal}
If
$(\pi , V)$ is a faithful representation of
$\bigoplus_{l=1}^r \h_{m,t^{d_l}}$ and
$d = \sum_{l = 1}^r d_l$,
 then
$$\dim V \geq m d + \left \lceil 2 \sqrt{d}\right\rceil.$$
In particular
$\displaystyle \mu\left(\h_{m,p}\right) \geq m  \deg p + \left \lceil 2 \sqrt{\deg p}\right\rceil$.
\end{theorem}

In what follows we shall assume that the exponents $d_1,\dots,d_r$
are fixed and we set $d = \sum_{l = 1}^r d_l$.

Let us define the following elements in $\bigoplus_{l=1}^r \h_{m,t^{d_l}}$:
\begin{align*}
X_{i,l}^j&=(0,\dots,0,X_i\otimes t^j,0,\dots,0), \\[1mm]
Y_{i,l}^j&=(0,\dots,0,Y_i\otimes t^j,0,\dots,0), \\[1mm]
Z_{l}^j&=(0,\dots,0,Z\otimes t^j,0,\dots,0),
\end{align*}
where the non-zero component of each element is in the $l^{\text th}$ coordinate.
Thus the set
\[
\Big\{X_{i,l}^j,\; Y_{i,l}^j,\; Z_{l}^j :
1 \leq i \leq m,\; 0 \leq j \leq d_l-1,\; 1 \leq l \leq r\Big\}
\]
is a  basis of $\bigoplus_{l=1}^r \h_{m,t^{d_l}}$.
The center $\z$  of $\bigoplus_{l=1}^r \h_{m,t^{d_l}}$ is spanned by the set
$\big\{Z_{l}^j :0 \leq j \leq d_l-1,\; 1 \leq l \leq r\big\}$.

\begin{lemma} \label{Z_r}
For any
$X \in \bigoplus_{l=1}^r \h_{m,t^{d_l}}$ such that $X\notin\z$  there exist
$Y \in \bigoplus_{l=1}^r \h_{m,t^{d_l}}$ such that $[X , Y] = Z_l^{d_l - 1}$
for some $l = 1 , \dots, r$.
\end{lemma}

\begin{proof}
Assume that
$$
X =\sum_{l=1}^r \sum_{j=0}^{d_l - 1}
\sum_{i = 1}^m a_{i,j,l} X_{i,l}^j +  b_{i,j,l} Y_{i,l}^j +  c_{j,l} Z_l^j
$$
for some $a_{i,j,l}, b_{i,j,l}, c_{j,l} \in \k$.
Since $X$ is not in the center of
$\bigoplus_{l=1}^r \h_{m,t^{d_l}}$,
then theree exists some $(i_0,j_0,l_0)$ such that either $a_{i_0,j_0,l_0}\ne0$ or $b_{i_0,j_0,l_0}\ne0$.
Assuming that $a_{i_0,j_0,l_0}\ne0$, let
$j_1= \min \{j : a_{i_0,j,l_0} \neq 0\}$ and let
$Y = \frac{1}{a_{i_0,j_1,l_0}} \,Y_{i_0,l_0}^{d_{l_0} - 1 -j_1}$ then
$[X , Y] = Z_{l_0}^{d_{l_0} - 1}$.
If $b_{i_0,j_0,l_0}\ne0$ the argument is analogous.
\end{proof}

\begin{lemma} \label{Lemma.dimSubalg}
Let $\g$ be a  Lie subalgebra of $\bigoplus_{l=1}^r
\h_{m,t^{d_l}}$ such that $Z_l^{d_l-1} \notin \g$ for all
$l = 1, \dots, r$. Then
$$\dim \g \leq m d + \dim \g \cap \mathfrak{z}.$$

\end{lemma}

\begin{proof}
Let
$\mathfrak{z}_0=\g \cap \mathfrak{z}$.
Since by hypothesis
$Z_l^{d_l - 1} \notin \mathfrak{z}_0$ for all
$l = 1, \dots, r$,
we may choose
a linear functional
$\alpha : \mathfrak{z} \rightarrow \k$ such that
$\alpha\mid_{\mathfrak{z}_0} = 0$  and
$\alpha(Z_l^{d_l-1})\neq 0$ for all $l = 1, \dots, r$.

Let
$\g_0$ be a complementary subspace of
$\mathfrak{z}_0$ in $\g$
let
$\tilde{\mathfrak{z}}$ be a complementary subspace of
$\mathfrak{z}_0$ in $\mathfrak{z}$,
and let
$\tilde{\g}$ be a complementary subspace of
$\g \oplus \tilde{\mathfrak{z}}$ in
$\bigoplus_{l=1}^r \h_{m,d_l}$.
Thus
$
\g = \g_0 \oplus \mathfrak{z}_0
$,
$
\z = \z_0 \oplus \tilde{\mathfrak{z}}
$
and
\begin{equation}\label{eq.decoposition}
\bigoplus_{l=1}^r \h_{m,t^{d_l}} =
\tilde{\g} \oplus \g_0 \oplus \mathfrak{z}_0 \oplus \tilde{\mathfrak{z}}
\end{equation}
as vector spaces.
Let $V = \tilde{\g} \oplus \g_0 $ and let
\begin{align*}
  B: V \times V  &\to        \k \\
     (X , Y)     &\mapsto    \alpha([X,Y])\text{ }.
\end{align*}
It is clear that
$B$ is an skew-symmetric  bilinear form on
$V$.

Let us prove that
$B$ is nondegenerate.
Given $X \in V$, $X\ne0$, we know from Lemma \ref{Z_r}
that there exist $Y \in \bigoplus_{l=1}^r \h_{m,t^{d_l}}$ such that
$[X , Y] = Z_l^{d_l - 1}$ for some
$l$.
If
$\tilde{Y} \in V$
is the projection of $Y$ to $V$ with respect to the
decomposition \eqref{eq.decoposition}, then
${B}(X , \tilde{Y}) \neq 0$.

Let us see now that
$\g_0$ is a ${B}$-isotropic subspace.
If
$X,Y \in \g_0$ then,
since
$\g$ is a Lie subalgebra, it follows that
$[X , Y] \in \mathfrak{z}_0$, and since
$\alpha\mid_{\mathfrak{z}_0} = 0$  we obtain that
$B(X , Y) = 0$.

Since
${B}$ is a nondegenerate  bilinear form on $V$ and
$\g_0$ is a ${B}$-isotropic subspace of
$V$, it follows that
$\dim \g_0 \leq \frac{\dim V}{2} = m d$ and therefore
$\dim \g \leq m d + \dim \z_0$.
\end{proof}

In order to prove Theorem \ref{TeoPrincipal}  we
need the following result
that gives some precise information about the structure of a
commuting family of nilpotent operators on a vector space.

\begin{theorem} \label{corleandro}
Let
$V$ be a finite dimensional vector space  and let
$\mathcal{N}$ be a non-zero abelian subspace of $\End (V)$
consisting of  nilpotent operators.
Then there exist a linearly independent set
$B=\{v_1,\dots,v_s\}\subset V$ and a decomposition
$\mathcal{N}=\mathcal{N}_1\oplus\dots\oplus\mathcal{N}_s$, with $\mathcal{N}_i\ne0$ for all $i$,
such that the applications
$F_i:\mathcal{N} \rightarrow V $ defined by
$F_i(N) = N(v_i)$ satisfy
\begin{enumerate}[\rm(1)]
\item  $F_i|_{\mathcal{N}_i}$ is injective for all
       $i=1\dots s$;
\item  $\mathcal{N}_{j}\subset\ker F_i$ for all
       $1\le i<j\le s$;
\item  $\mathcal{N}_{j}V \subset \im F_i|_{\mathcal{N}_i}$ for all
       $1\le i<j\le s$.
\end{enumerate}
Furthermore, given a finite subset
$\{N_1,\dots,N_q\}$ of non-zero operators in $\mathcal{N}$, the vector $v_1$
can be chosen so that $N_k(v_1) \neq 0$ for all $k=1,\dots, q$.
\end{theorem}
We first prove the following lemma.
\begin{lemma} \label{lemav1}
Let
$V$ be a finite dimensional vector space over
$\k$,  let
$\mathcal{F}$ a non-zero subspace of $\End(V)$
and let $r= \max \{ \dim \mathcal{F}v : v \in V\}$.
Then for any finite subset $\{T_1, T_2,\dots,T_q\} \subseteq \mathcal{F}$,
such that $T_i\ne0$ for all $i=1,\dots,q$, there exist
$v_0 \in V$ such that
$r = \dim \mathcal{F}v_0 \text{ and } T_i(v_0) \neq 0$ for all
$i = 1, \dots, q$.
\end{lemma}

\begin{proof}
We will prove the lemma by induction on $q$.
It is clear that the lemma is true in the case
$q = 0$. Let
$\{T_1, \dots,T_q, T_{q+1}\} \subseteq \mathcal{F}$, by inductive hypothesis there exist
$v'_0 \in V$ such that
$r = \dim \mathcal{F}v'_0$ and
$T_i(v'_0) \neq 0$ for all
$i = 1, \dots, q$. If $T_{q+1}(v'_0) \neq 0$ we take $v_0= v'_0$.

Suppose that $T_{q+1}(v'_0) = 0$. Since
$T_{q+1} \neq 0$ there exist
$w \in V$ such that
$T_{q+1}(w) \neq 0$. Let us take
$\tilde{T}_1, \dots, \tilde{T}_r\in\mathcal{F}$ such that
$\big\{\tilde{T}_1(v'_0) , \dots, \tilde{T}_r(v'_0)\big\}$ is a basis of
$\mathcal{F}v'_0$.

We now claim that there exists $t_0$ such that

\begin{enumerate}[(1)]
\item the set $\{\tilde{T}_1(v'_0 + t_0w), \dots, \tilde{T}_r(v'_0 + t_0w)\}$
is linearly independent, and

\item $T_i(v'_0 + t_0w) \neq 0$ for all
$i = 1, \dots, q+1$.
\end{enumerate}
In order to prove this fact, let $B$ be a basis of $V$.
Let $A_t$ be the matrix whose columns are the coordinates of the
vectors $\tilde{T}_1(v'_0 + tw), \dots, \tilde{T}_r(v'_0 + tw)$ and
let
$a(t)$ be the $r\times r$ minor of $A_t$ such that
$a(0)\ne0$
(since $\big\{\tilde{T}_1(v'_0) , \dots, \tilde{T}_r(v'_0)\big\}$ is linearly independent
the existence of this minor is granted).
For $i = 1, \dots, q$,
let $p_i(t)$ be a coordinate of ${T}_i(v'_0 + tw)$ such that
$p_i(0)\ne0$ and
let $p_{q+1}(t)$ be a coordinate of ${T}_{q+1}(v'_0 + tw)$ such that
$p_{q+1}(1)\ne0$ (recall that we assumed that ${T}_{q+1}(v'_0 )=0$).
Now $\{a(t),p_1(t),\dots,p_{q+1}(t)\}$ is a finite set of
non-zero polynomials and thus there exist $t_0$ such that
non of them vanish at $t_0$.
For this $t_0$ conditions (1) and (2) are verified and
taking $v_0=v'_0 + t_0w$ we complete the inductive argument.
\end{proof}

\begin{proof}[Proof of Theorem \ref{corleandro}]
We shall proceed by induction on
$\dim \mathcal{N}$.
If
$\dim \mathcal{N} = 1$, let $N_0$ be any non-zero operator in $\mathcal{N}$ and
let $v_1$ be a vector such that $N_0(v_1)\ne0$.
If we take $B=\{v_1\}$ and
$\mathcal{N}_1=\mathcal{N}$ then
$F_1|_{\mathcal{N}_1}$ is injective
and conditions (2) and (3) are empty.
It is clear that if $N_1,\dots,N_q$ are non-zero operators in
$\mathcal{N}$ then $N_k(v_1) \neq 0$ for all $k=1,\dots, q$.

Now assume that the theorem is true for any non-zero abelian subspace of
$\End (V)$ of dimension less than $\dim \mathcal{N}$.

Let
$r = \max\{\dim \mathcal{N}v : v \in V\}>0$.
By Lemma \ref{lemav1}, there exist
$v_1 \in V$ such that $r = \dim \mathcal{N}v_1$ and
$N_k(v_1) \neq 0$ for all $k=1,\dots, q$.
Let
$F_1:\mathcal{N} \rightarrow  V$ defined by
$F_1(N)= N(v_1)$.
If $F_1$ is injective then
 we take $B=\{v_1\}$ and
$\mathcal{N}_1=\mathcal{N}$.
We obtain that $F_1|_{\mathcal{N}_1}$ is injective
and since conditions (2) and (3) are empty, we are done.

Otherwise, let
$\mathcal{N}'= \ker F_1$.
Since
$r = \dim \mathcal{N}v_1>0$, we have
$\dim\mathcal{N}'<\dim\mathcal{N}$.
By the inductive hypothesis, there exist a linearly independent set
$B'=\{v_2, \dots, v_s\} \subset V$ and a decomposition
$\mathcal{N}'=\mathcal{N}_2\oplus\dots\oplus\mathcal{N}_s$,
with $\mathcal{N}_i\ne0$ for $i=2,\dots,s$, such that
\begin{enumerate}[\rm($1'$)]
  \item  $F_i|_{\mathcal{N}_i}$ is injective for all
         $i=2\dots s$;
  \item  $\mathcal{N}_{j} \subset \ker F_i$ for all
         $2\le i<j\le s$;
  \item  $\mathcal{N}_{j}V \subset \im F_i|_{\mathcal{N}_i}$ for all
         $2\le i<j\le s$.
\end{enumerate}
Let us prove that $B = \{v_1, v_2, \dots, v_s\} $ is
a linearly independent set.
Since $B'$ is linearly independent we must show that
$v_1 \not\in \k B'$.
If
$v_1 \in \k B'$ then there exist
$a_j \in \k$, not all of them zero, such that
$v_1 = \sum_{j=2}^s a_j v_j$. Let
$j_0 = \max \{j : a_j \neq 0 \}$ and
let $N \in \mathcal{N}_{j_0}\subset\mathcal{N}' $ be a non-zero operator.
If we apply
$N$ to both sides of
$v_1 = \sum_{j=2}^s a_j v_j$ we obtain zero on the
left hand side (since $N \in \mathcal{N}' $), but
from conditions $(1')$ and $(2')$,
 we obtain $ a_{j_0} N(v_{j_0})\ne0$ on the right hand side,
 which is a contradiction.

 Let
$\mathcal{N}_1$ be a direct complement of
$\mathcal{N}'$ in $\mathcal{N}$.
It is clear that
$\mathcal{N}=\mathcal{N}_1\oplus\dots\oplus\mathcal{N}_s$, $\mathcal{N}_i\ne0$ for all $i$,
and that conditions (1) and (2) are verified.
In order to finish the inductive step we must prove that
conditions (3) is verified.
In fact we only need to show that
$N'v'\in\im F_1|_{\mathcal{N}_1}$ for all
$N'\in\mathcal{N}'$ and $v'\in V$.

Since
$\dim \mathcal{N} = \dim \ker F_1 + \dim \im F_1 = \dim \mathcal{N}' + r$,
we know that
$\dim \mathcal{N}_1 = r$.
Let
$\{\tilde{N}_1,\dots,\tilde{N}_r\}$ be a basis of
${\mathcal{N}}_{1}$.
Given arbitrary $N'\in\mathcal{N}'$ and $v'\in V$
we must show that
$N'v'\in\k\big\{\widetilde{N}_1v_1,\dots,\widetilde{N}_rv_1\big\}$
which in turns is equivalent to prove that
$\big\{N'v',\widetilde{N}_1v_1,\dots,\widetilde{N}_rv_1\big\}$ is linearly dependent.
By the definition of
$v_1$, the set
$\big\{N'(v_1+tv'),\tilde{N}_1(v_1+tv'),\dots,\tilde{N}_r(v_1+tv')\big\}$
is linearly dependent for all
$t \in \k$. Since
$N'v_1=0$, we have that
\[
B' = \big\{N'v',\widetilde{N}_1(v_1+tv'),\dots,\widetilde{N}_r(v_1+tv')\big\}
\]
is linearly dependent for all
$t \ne0$, and therefore it is linearly dependent for $t=0$,
as we wanted to prove.
\end{proof}

We are now ready to prove the main result of this section.

\begin{proof}[Proof of Theorem \ref{TeoPrincipal}.]
Let
$(\pi , V)$ be a faithful representation of
$\bigoplus_{l=1}^k \h_{m,t^{d_l}}$.
By Theorem \ref{thm.nilrep} we may assume that $\pi$ is a nilrepresentation.

We apply Theorem \ref{corleandro} to the
subspace $\mathcal{N}=\pi(\z)$.
We obtain
a linearly independent set $B=\{v_1,\dots,v_s\}\subset V$
and decomposition
$\mathcal{N}=\mathcal{N}_1\oplus\dots\oplus\mathcal{N}_s$, $\mathcal{N}_i\ne0$ for all $i$,
such that the applications
$F_i:\mathcal{N} \rightarrow V $ defined by
$F_i(N) = N(v_i)$ satisfy
\begin{enumerate}[\rm(1)]
\item  $F_i|_{\mathcal{N}_i}$ is injective for all
       $i=1\dots s$;
\item  $\mathcal{N}_{j}\subset\ker F_i$ for all
       $1\le i<j\le s$;
\item  $\mathcal{N}_{j}V \subset \im F_i|_{\mathcal{N}_i}$ for all
       $1\le i<j\le s$.
\end{enumerate}
We additionally require that $\pi(Z_l^{d_l-1})v_1\ne0$
for all $l=1,\dots, r$.

Let $\phi$ be the linear map
\[
  \phi:\bigoplus_{l=1}^r \h_{m,t^l}   \to        V,\quad
          \phi(X)=\pi(X)v_1.
\]
Note that $\phi|_{\z}=F_1\circ\pi$.
We claim that
\begin{enumerate}[\rm(i)]
\item  $\dim \im\phi + \dim \ker F_1 \geq (m+1)d$.
\item  $\im \phi\cap\k B=0$, and thus $\dim V\ge s+\dim \im \phi$.
\item  $d\leq s\, \dim\im F_1$ and thus
$\left \lceil 2 \sqrt{d} \right\rceil \leq
s + \dim\im F_1$.

\end{enumerate}

\noindent
\textit{Proof of (i)}. It is clear that $\ker\phi$
is a subalgebra of $\bigoplus_{l=1}^k \h_{m,t^l}$ such that
$Z_l^{d_l-1} \notin \ker\phi$.
Since $\pi\big((\ker\phi)\cap\z\big)=\ker F_1$, we obtain
from Lemma \ref{Lemma.dimSubalg}
that
\[
\dim \ker\phi \leq m d + \dim\ker F_1.
\]
Since $\dim\ker\phi+\dim\im\phi=(2m+1)d$, we obtain part (i).

\medskip
\noindent
\textit{Proof of (ii)}.
Let $v \in \im \phi\cap\k B$.
Since
$v \in \im \phi$ there exists
$X \in \bigoplus_{l=1}^r \h_{m,t^l}$ such that
$\pi(X)(v_1) = v$
and thus there exist
$a_1,\dots,a_s \in \k$  such that
\begin{equation}\label{eq.imagen}
\pi(X)v_1 = \sum_{i = 1}^s a_i v_i.
\end{equation}
We must prove that $a_i=0$ for all $i$.
Assume that $a_i\ne0$ for some $i$ and  let $i_0 = \max \{i : a_i \neq 0 \}$.
Since $\pi(X)$ is a nilpotent endomorphism on $V$, its only eigenvalue is zero and thus $i_0 > 1$.
Now let $N \in \mathcal{N}_{i_0}$ be a non-zero operator and let us  apply
$N$ to both sides of equation \eqref{eq.imagen}. We obtain zero on the
left hand side, since
$N \in \pi(\mathfrak{z})$ and
$i_0 >1$. But from conditions $(1)$ and $(2)$, we obtain
$a_{i_0} N(v_{i_0})\ne0$ on the right hand side, which is a contradiction.

\medskip
\noindent
\textit{Proof of (iii)}.
Part (1) and (3) combined imply that
$\dim\mathcal{N}_x\ge\dim\mathcal{N}_y$ if
$x<y$.
In particular $\dim\mathcal{N}_1\ge\dim\mathcal{N}_j$ for all $j=1,\dots,s$ and thus
\[
d=\dim\mathcal{N}=\sum_{j=1}^s\dim\mathcal{N}_j\le s\,\dim\mathcal{N}_1=s\,\dim\im F_1.
\]
Since
$
\min\{a+b:a,b\in \mathbb{N} \text{ and } ab\ge d\}= \left \lceil 2\sqrt{d} \right\rceil
$
 for all $d\in \mathbb{N}$,
 we obtain part (iii).

From part (i) and (ii) it follows that
\begin{equation*}
\dim V + \dim \ker F_1 \geq (m+1)d + s,
\end{equation*}
and combining it with part (iii)
we obtain
\begin{equation*}
\dim V + \dim \ker F_1 + \dim\im F_1\geq (m+1)d + \left \lceil 2 \sqrt{d} \right\rceil.
\end{equation*}
Finally, since
$\dim \ker F_1 + \dim\im F_1=d$ we obtain
\begin{equation*}
\dim V  \geq md + \left \lceil 2 \sqrt{d} \right\rceil
\end{equation*}
as we wanted to prove.
\end{proof}

We close this section with the following corollary that is
equivalent to the Schur's Theorem mentioned in the introduction.

\begin{corollary}\label{Coro.Schur}
Let
$V$ be a finite dimensional vector space and let
$\mathcal{N}$ be a non-zero abelian subspace of $\End (V)$
consisting of  nilpotent operators.
Then
\[
\dim V\ge\left\lceil2\sqrt{\dim\mathcal{N}}\right\rceil.
\]
\end{corollary}
\begin{proof}
  Let $B$, $\mathcal{N}_1$, $s$ and $F_1$ as in Theorem \ref{corleandro}.
  Then, by the same argument used in the proof of Theorem \ref{TeoPrincipal}
  in it items (ii) and (iii)
  we obtain that
\medskip
\begin{enumerate}
\item[(ii)'] $\im F_1\cap\k B=0$, and thus $\dim V\ge s+\dim \im F_1$, and

\item[(iii)'] $\dim\mathcal{N}\leq s\, \dim\im F_1$ and thus
$\left \lceil 2 \sqrt{d} \right\rceil \leq
s + \dim\im F_1$.
\end{enumerate}
Therefore  $\dim V\ge\left\lceil2\sqrt{\dim\mathcal{N}}\right\rceil$.
\end{proof}

\addcontentsline{toc}{chapter}{Bibliograf{\'\i}a}






\end{document}